# Using a DPG method to validate DMA experimental calibration of viscoelastic materials

Federico Fuentes[1], Leszek Demkowicz[1], and Aleta Wilder[2]

[1]The Institute for Computational Engineering and Sciences (ICES), The University of Texas at Austin, 201 E 24th St, Austin, TX 78712, USA

[2]Center for Energy and Environmental Resources (CEER), The University of Texas at Austin, 10100 Burnet Road, Bldg. 133, Austin, TX 78758, USA

**Abstract**

A discontinuous Petrov-Galerkin (DPG) method is used to solve the time-harmonic equations of linear viscoelasticity. It is based on a "broken" primal variational formulation, which is very similar to the classical primal variational formulation used in Galerkin methods, but has additional "interface" variables at the boundaries of the mesh elements. Both the classical and broken formulations are proved to be well-posed in the infinite-dimensional setting, and the resulting discretization is proved to be stable. A full $hp$-convergence analysis is also included, and the analysis is verified using computational simulations. The method is particularly useful as it carries its own natural arbitrary-$p$ a posteriori error estimator, which is fundamental for solving problems with localized solution features. This proves to be useful when validating calibration models of dynamic mechanical analysis (DMA) experiments. Indeed, different DMA experiments of epoxy and silicone resins were successfully validated to within 5% of the quantity of interest using the numerical method.

## 1 Introduction

The advent of minimum residual methods to solve general well-posed linear variational formulations has made it possible to revisit various problems in the literature, especially those with notable numerical stability issues, such as convection-diffusion [17, 11, 23, 9, 53]. The discontinuous version of these methods, called the discontinuous Petrov-Galerkin (DPG) methodology, uses broken (i.e. discontinuous) test spaces and extra interface variables to optimize a discrete Riesz map inversion in order to approximate a solution minimizing a variational residual and aiming to reproduce the best possible stability. Despite coming at a relatively high computational cost, the DPG methodology not only is very general and leads to stable and convergent numerical methods, but has other numerous advantages, such as always resulting in positive-definite stiffness matrices; a natural a posteriori error estimator that can be used to drive adaptive schemes allowing problems with localized solution features to be tackled; and an accessible way of coupling with other numerical methods across boundaries. For example, the equations of static linear elasticity have been flexibly solved with the methodology using several different variational formulations [43], which were then coupled to capitalize on the fact that some formulations were robust in the incompressible limit while others were computationally more efficient [27]. Similarly, different features of the methodology have been taken advantage of while solving problems related to fluid flow [59, 8, 26, 44], wave propagation [69, 31, 18], electromagnetism [6], elasticity [43, 7, 5], transmission problems in unbounded domains [38, 29], and even optical fibers [19], among others.

The aim of this article is twofold. First, it is to implement a primal formulation of time-harmonic linear viscoelasticity with a DPG method. Second, it is to use such an implementation to validate calibration data directly from dynamic mechanical analysis (DMA) experimental results of the dynamic Young's modulus of two different thermoset resins. Many problems in viscoelasticity have local solution features in the stress or displacement, and the a posteriori error estimate is a very useful trait of the general DPG methodology which can be exploited in those cases. In fact, such local solution features will be observed when simulating the experimental results.

The article is outlined as follows. In Section 2 the equations of viscoelasticity along with the relevant variational formulations are introduced and proved to be well-posed. In Section 3 minimum residual methods and the DPG methodology are described and the corresponding discrete numerical method is shown to be stable and $hp$-convergent. In Section 4, numerical results that verify the numerical scheme are presented along with a validation of calibration models from DMA experimental results.

## 2 Primal variational formulations for viscoelasticity

### 2.1 Linear viscoelasticity

The classical linear viscoelasticity equations are solved in this work. The constitutive model originally was developed by Boltzmann [4] and Volterra [67], but later recast more rigorously as a linearization of nonlinear continuum mechanics under the additional assumption of a dependence of the stress on the deformation history [13, 12, 54, 55]. In the time domain, the first-order system describing a viscoelastic material with constant density $\rho > 0$ in a domain $\Omega \subseteq \mathbb{R}^3$ is

$$\begin{cases} \rho \ddot{u} = \operatorname{div} \sigma + f \,, \\ \sigma = \dot{\mathsf{C}} * \varepsilon = \int_{-\infty}^{\infty} \dot{\mathsf{C}}(s) \colon \varepsilon(\,\cdot\, - s) \, \mathrm{d}s \,, \end{cases} \quad (2.1)$$

where the displacement $u$ and stress $\sigma$ are unknown, $f$ is a known body force, and the engineering strain is defined in terms of $u$ as $\varepsilon = \frac{1}{2}(\nabla u + \nabla u^{\mathsf{T}})$. Meanwhile, the viscoelastic stiffness tensor $\mathsf{C}$ is in general not only a function in space, but also in time. With the typical assumption of $\mathsf{C}(t) = 0$ for times $t < 0$, this leads to the distributional derivative $\dot{\mathsf{C}}(t) = \mathsf{C}(0)\delta_0(t) + \dot{\mathsf{C}}^+(t)H_0(t)$, where $\dot{\mathsf{C}}^+ = \frac{\mathrm{d}\mathsf{C}|_{(0,\infty)}}{\mathrm{d}t}$, $\delta_0$ is the Dirac delta, and $H_0$ is the Heaviside step function. This leads to the expression $\sigma = \mathsf{C}(0) \colon \varepsilon + \int_0^{\infty} \dot{\mathsf{C}}(s) \colon \varepsilon(\,\cdot\, - s) \, \mathrm{d}s$, which is commonly found in the literature [36, 13, 25]. The relaxation or equilibrium stiffness tensor is $\mathsf{C}^{\infty} = \lim_{t \to \infty} \mathsf{C}(t)$. The classical case of linear elasticity occurs when $\mathsf{C}(t) = \mathsf{C}^{\infty} H_0(t)$ leading to $\sigma = \mathsf{C}^{\infty} \colon \varepsilon$.

In practice, many applications occur in a vibrating environment, so considering the time-harmonic case is natural. This also has the advantage of avoiding the computation of any convolutions, since $\dot{\mathsf{C}} * \varepsilon$ becomes a product after using the Fourier transform. As usual, the stiffness tensor is assumed to have minor and major symmetries, so that $\mathsf{C}_{ijkl} = \mathsf{C}_{ijlk} = \mathsf{C}_{jikl} = \mathsf{C}_{klij}$ and as a result $\mathsf{C}_{ijkl}\tau_{kl} = \mathsf{C}_{ijkl}\frac{1}{2}(\tau_{kl} + \tau_{lk})$ for any second-order tensor $\tau$. In particular $\dot{\mathsf{C}} * \varepsilon = \dot{\mathsf{C}} * \nabla u$. Thus, substituting the constitutive model for the stress into the conservation of momentum, and considering the time-harmonic case at angular frequency $\omega$, yields the second-order equation,

$$-\omega^2 \rho u - \operatorname{div}(\mathsf{C}^* \colon \nabla u) = f \,, \quad (2.2)$$

where, given $x \in \Omega$, the complex-valued $u(x, \omega)$, $f(x, \omega)$ and $\mathsf{C}(x, \omega)$ are the corresponding Fourier transforms of the time-dependent displacement, force and stiffness tensor; and where $\mathsf{C}^*(x, \omega) = \mathrm{i}\omega \mathsf{C}(x, \omega)$ is defined as the dynamic stiffness tensor. Note that in the limiting case of linear elasticity, $\mathsf{C}^* = \mathsf{C}^\infty$, so $\mathsf{C}^*$ is no longer complex-valued or $\omega$-dependent. For isotropic materials, the stiffness tensor explicitly takes the form $\mathsf{C}_{ijkl} = \lambda \delta_{ij}\delta_{kl} + \mu(\delta_{ik}\delta_{jl} + \delta_{il}\delta_{jk})$ and similarly with $\mathsf{C}^*$ in terms of $\lambda^* = i\omega\lambda$ and $\mu^* = i\omega\mu$, where in the latter expression $\lambda$ and $\mu$ are the Fourier transforms of the time-dependent Lamé parameters. Moreover, $G^* = \mu^*$ is the dynamic shear modulus, $K^* = \lambda^* + \frac{2}{3}\mu^*$ is the dynamic bulk modulus, $E^* = \frac{\mu^*(3\lambda^* + 2\mu^*)}{\lambda^* + \mu^*}$ is the dynamic Young's modulus, and $\nu^* = \frac{\lambda^*}{2(\lambda^* + \mu^*)}$ is the dynamic Poisson's ratio. Notably, $E^*$ is a nonlinear function of $\lambda^*$ and $\mu^*$, implying that in general it is not the Fourier transform of $\frac{d}{dt}\frac{\mu(t)(3\lambda(t) + 2\mu(t))}{\lambda(t) + \mu(t)}$. A similar assertion follows for $\nu^*$. Thus, one should be careful when speaking of the time-dependent Young's modulus and Poisson's ratio in three-dimensional viscoelasticity as even different definitions derived from physical principles exist in the literature [46, §5.7].

The goal is to solve the second-order equation in (2.2) for the unknown displacement, provided the forcing and the dynamic stiffness tensor of the material are known throughout the domain $\Omega \subseteq \mathbb{R}^3$ at the angular frequency $\omega$. For this to be possible, boundary conditions need to be specified, so it will be assumed that the boundary is partitioned into relatively open subsets $\Gamma_u$ and $\Gamma_\sigma$ satisfying $\overline{\Gamma_u \cup \Gamma_\sigma} = \partial\Omega$ and $\Gamma_u \cap \Gamma_\sigma = \varnothing$, where displacement and traction boundary conditions are set by the known functions $u = u^{\Gamma_u}$ and $(\mathsf{C}^* : \nabla u) \cdot \hat{n} = \sigma_n^{\Gamma_\sigma}$ on $\Gamma_u$ and $\Gamma_\sigma$ respectively, with $\hat{n}$ being the outward normal at $\partial\Omega$. From now on it will be assumed that $\Gamma_u \neq \varnothing$ and $\Omega$ is bounded and Lipschitz.

## 2.2 Classical primal formulation

The usual approach to solve the second-order equation is to multiply by a smooth enough test function that vanishes at $\Gamma_u$, and then integrate by parts once. This yields the sesquilinear and conjugate-linear forms,

$$b_0(u,v) = -\omega^2 \rho(u,v)_\Omega + (\mathsf{C}^* : \nabla u, \nabla v)_\Omega = (f,v)_\Omega + \langle \sigma_n^{\Gamma_\sigma}, v \rangle_{\partial\Omega} = \ell(v), \qquad (2.3)$$

for all test functions $v$. Here $(w_1, w_2)_K = \int_K \text{tr}(w_2^\dagger w_1)\, dK$ for $K \subseteq \Omega$, where tr is the usual trace of a matrix and $w_2^\dagger$ is the conjugate transpose of $w_2$, so that depending on whether $w_1$ and $w_2$ take scalar, vector or matrix values, $\text{tr}(w_2^\dagger w_1)$ will be $w_1 \bar{w}_2$, $w_1 \cdot \bar{w}_2$ or $w_1 : \bar{w}_2$ respectively. Similarly, if $\sigma_n^{\Gamma_\sigma}$ and $v$ are smooth enough, $\langle \sigma_n^{\Gamma_\sigma}, v \rangle_{\partial K}$ would be a boundary integral over $\partial K$ for $K \subseteq \Omega$ (where, abusing notation, $\sigma_n^{\Gamma_\sigma}$ is understood as any extension from $\Gamma_\sigma$ to $\partial K$).

To prove the convergence and stability of any numerical method aiming to solve (2.3), usually determining well-posedness of the underlying non-discrete equations is either necessary or extremely useful. For this, a deeper understanding of the functional spaces used as trial and test spaces is required. Indeed, when $u^{\Gamma_u} = 0$, the natural choice of space for $u$ and $v$ is $H^1_{\Gamma_u}(\Omega) = \{u \in H^1(\Omega) \mid u|_{\Gamma_u} = 0\}$, where for any $K \subseteq \Omega$,

$$\begin{aligned} L^2(K) &= \{u : K \to \mathbb{C}^3 \mid \|u\|^2_{L^2(K)} = (u,u)_K < \infty\}, \\ H^1(K) &= \{u : K \to \mathbb{C}^3 \mid \|u\|^2_{H^1(K)} = (u,u)_K + (\nabla u, \nabla u)_K < \infty\}. \end{aligned} \qquad (2.4)$$

When, $u^{\Gamma_u} \neq 0$, the final displacement takes the form $u_f = u + \widetilde{u}^{\Gamma_u}$, where $u \in H^1_{\Gamma_u}(\Omega)$ and $\widetilde{u}^{\Gamma_u} \in H^1(\Omega)$ is an extension of $u^{\Gamma_u}$ to $\Omega$. For simplicity consider $u^{\Gamma_u} = 0$, let $U = H^1_{\Gamma_u}(\Omega)$ and assume that for some $C > 0$, $|\ell(v)| \leq C\|v\|_{H^1(K)}$ for all $v \in U$, so that $\ell \in U'$, with $U'$ being the space of *conjugate*-linear continuous functionals with domain $U$. Then, solving (2.3) is equivalent to the problem

$$\begin{cases} \text{Find } u \in U, \\ b_0(u,v) = \ell(v), \quad \text{for all } v \in U, \end{cases} \qquad (2.5)$$

and the goal is to prove this equation is well-posed in the sense of Hadamard, so that there is a guaranteed existence of a unique solution depending continuously upon the forcing and boundary conditions (encoded in $\ell$). The proof is presented in what remains of the section, where a bounded $\Omega \subseteq \mathbb{R}^3$ and $\Gamma_u \neq \varnothing$ are assumed throughout. It is based on the use of the Fredholm alternative and the theory of Gelfand triples in the same spirit as [37, 56, 41].

**Lemma 2.1.** *Let $b_{\mathsf{C}^*}(u,v) = (\mathsf{C}^* : \nabla u, \nabla v)_\Omega$, for $u, v \in H^1_{\Gamma_u}(\Omega)$ and with $\mathsf{C}^*$ being a fourth-order tensor with major and minor symmetries satisfying $|\bar{\varepsilon} : \mathfrak{Re}(\mathsf{C}^*) : \varepsilon| > 0$ for all symmetric second-order tensors $\varepsilon \neq 0$. Then, $|b_{\mathsf{C}^*}(u,u)| \geq \alpha \|u\|^2_{H^1(\Omega)}$ for some $\alpha > 0$.*

*Proof.* First note that $\overline{\nabla u} : \mathsf{C}^* : \nabla u = \bar{\varepsilon} : \mathfrak{Re}(\mathsf{C}^*) : \varepsilon + i\bar{\varepsilon} : \mathfrak{Im}(\mathsf{C}^*) : \varepsilon$ with $\varepsilon = \frac{1}{2}(\nabla u + \nabla u^\mathsf{T})$. The major symmetry of $\mathsf{C}^*$ clearly implies that both $\bar{\varepsilon} : \mathfrak{Re}(\mathsf{C}^*) : \varepsilon$ and $\bar{\varepsilon} : \mathfrak{Im}(\mathsf{C}^*) : \varepsilon$ are real-valued, so that

$$|b_{\mathsf{C}^*}(u,u)|^2 = |(\mathfrak{Re}(\mathsf{C}^*) : \varepsilon, \varepsilon)_\Omega|^2 + |(\mathfrak{Im}(\mathsf{C}^*) : \varepsilon, \varepsilon)_\Omega|^2 \geq |(\mathfrak{Re}(\mathsf{C}^*) : \varepsilon, \varepsilon)_\Omega|^2.$$

Due to the symmetries, $\mathfrak{Re}(\mathsf{C}^*)$ and $\varepsilon$ may be reinterpreted in Voigt notation as a symmetric $6 \times 6$ matrix and a vector in $\mathbb{C}^6$ respectively, so that the Rayleigh quotient of the Voigt-matrix $\mathfrak{Re}(\mathsf{C}^*)$ takes the form $\bar{\varepsilon} : \mathfrak{Re}(\mathsf{C}^*) : \varepsilon / (\bar{\varepsilon} : \varepsilon + 2|\varepsilon_{12}|^2 + 2|\varepsilon_{13}|^2 + 2|\varepsilon_{23}|^2)$. By hypothesis, 0 is not in its range, implying that its range is fully positive or fully negative, and that $|(\mathfrak{Re}(\mathsf{C}^*) : \varepsilon, \varepsilon)_\Omega| = \int_\Omega |\bar{\varepsilon} : \mathfrak{Re}(\mathsf{C}^*) : \varepsilon| \, d\Omega$. If the range is positive, the Rayleigh quotient yields $|\bar{\varepsilon} : \mathfrak{Re}(\mathsf{C}^*) : \varepsilon| \geq \lambda_{\min}(\bar{\varepsilon} : \varepsilon + 2|\varepsilon_{12}|^2 + 2|\varepsilon_{13}|^2 + 2|\varepsilon_{23}|^2) \geq \lambda_{\min}\bar{\varepsilon} : \varepsilon$, where $\lambda_{\min} > 0$ is the smallest eigenvalue of the Voigt-matrix $\mathfrak{Re}(\mathsf{C}^*)$. Similarly if the range is negative, so that in any case $|b_{\mathsf{C}^*}(u,u)| \geq \alpha(\varepsilon,\varepsilon)_\Omega$ for some $\alpha > 0$. The result follows because Korn's and Poincaré inequalities ($\Gamma_u \neq \varnothing$) imply that for all $u \in H^1_{\Gamma_u}(\Omega)$, $(\varepsilon,\varepsilon)_\Omega \geq \alpha\|u\|^2_{H^1(\Omega)}$ for some $\alpha > 0$. □

**Remark 2.1.** *In the case of isotropic materials, the conditions on $\mathsf{C}^*$ are equivalent to $\mathfrak{Re}(G^*)\mathfrak{Re}(K^*) > 0$. The physically relevant case is when both the storage shear and bulk moduli are positive, $\mathfrak{Re}(G^*) > 0$ and $\mathfrak{Re}(K^*) > 0$, but exotic exceptions do exist where the storage bulk modulus may be negative [47]. Curiously, if $\mathfrak{Re}(G^*) \neq 0$ and $\mathfrak{Re}(K^*) = 0$, then $I : \mathfrak{Re}(\mathsf{C}^*) : I = 0$, but the coercive inequality $|b_{\mathsf{C}^*}(u,u)| \geq \alpha\|u\|^2_{H^1(\Omega)}$ still holds, because $|(\mathfrak{Re}(\mathsf{C}^*) : \varepsilon, \varepsilon)_\Omega| = 2\mathfrak{Re}(G^*)(\varepsilon^D, \varepsilon^D)_\Omega$, where $\varepsilon^D = \varepsilon - \frac{1}{3}\mathrm{tr}(\varepsilon)I$ is the deviatoric part of the strain. Then, all that remains is to apply a recently proved and more general version of Korn's inequality, $(\varepsilon^D, \varepsilon^D)_\Omega \geq \alpha(\nabla u, \nabla u)_\Omega$ for all $u \in H^1_{\Gamma_u}(\Omega)$ and $\alpha > 0$ [52].*

**Remark 2.2.** *In the particular case of static linear elasticity, $\mathsf{C}^* = \mathfrak{Re}(\mathsf{C}^*) = \mathsf{C}^\infty$ and the primal formulation is that of finding $u \in H^1_{\Gamma_u}(\Omega)$ such that $b_{\mathsf{C}^*}(u,v) = \ell(v)$ for all $v \in H^1_{\Gamma_u}(\Omega)$. Thus, a straightforward application of the Lax-Milgram theorem yields the well-posedness of the static linear elasticity equation provided $|\varepsilon : \mathsf{C}^\infty : \varepsilon| > 0$ for all symmetric strains $\varepsilon \neq 0$. If the material is isotropic this implies $G^\infty K^\infty > 0$, and in particular, the equations are well-posed for positive shear and bulk moduli. For even more general conditions (in terms of the compliance tensor, $\mathsf{S}^\infty = (\mathsf{C}^\infty)^{-1}$) under which static linear elasticity remains well-posed, see [1].*

**Theorem 2.1.** *Let $U = H^1_{\Gamma_u}(\Omega)$ and consider the problem of finding $u \in U$ such that $b_0(u,v) = \ell(v)$ for all $v \in U$, where $b_0(u,v) = -\omega^2\rho(u,v)_\Omega + (\mathsf{C}^* : \nabla u, \nabla v)_\Omega$ and $\ell \in U'$, and where $|\bar{\varepsilon} : \mathfrak{Re}(\mathsf{C}^*) : \varepsilon| > 0$ for all symmetric second-order tensors $\varepsilon \neq 0$. Then, for each value of $\omega$, either there exists $0 \neq u \in U$ such that $b_0(u,v) = 0$ for all $v \in U$, or, given any $\ell \in U'$, there exists a unique solution $u \in U$ solving $b_0(u,v) = \ell(v)$ for all $v \in U$ which satisfies $\|u\|_U \leq C\|\ell\|_{U'}$ for a $C > 0$ independent of the choice of $\ell$. Furthermore, the former case, where infinitely many solutions of the form $\beta u \in U$ for $\beta \in \mathbb{C}$ exist, only holds for a countable set of values of $\omega$ which has no accumulation points.*

*Proof.* First define the linear operator $B : U \to U'$ as $\langle Bu, v \rangle = b_{\mathsf{C}^*}(u,v) = (\mathsf{C}^* : \nabla u, \nabla v)_\Omega$ for all $v \in U$. Lemma 2.1 implies that $B$ is bounded below, $\|Bu\|_{U'} \geq \alpha \|u\|_U$, with some $\alpha > 0$, so that $B$ is injective and has closed range $\mathsf{R}(B) = \{\ell \in U' \mid \ell|_{U_{00}} = 0, U_{00} = \{v \in U \mid b_{\mathsf{C}^*}(u,v) = 0 \,\forall u \in U\}\}$. Again by Lemma 2.1, $U_{00} = \{0\}$ and $\mathsf{R}(B) = U'$, so the open mapping theorem implies $B^{-1} : U' \to U$ is bounded. Assume the embedding $\iota : U \to U'$, defined naturally as $\langle \iota u, v \rangle = (u,v)_\Omega$ for all $v \in U$, is compact, so that the operator $K = \iota B^{-1} : U' \to U'$ is a compact operator, with range $\mathsf{R}(K) = \iota(U)$. Given $\omega \neq 0$ (so $\omega^2 \rho \neq 0$), the Fredholm alternative applies to $K - \frac{1}{\omega^2 \rho} \mathrm{id}$. Therefore, either there exists $0 \neq \iota u = v \in \mathsf{R}(K)$ such that $Kv - \frac{1}{\omega^2 \rho} v = 0$ and so $-\omega^2 \rho B \iota^{-1}(Kv - \frac{1}{\omega^2 \rho} v) = -\omega^2 \rho \iota u + Bu = 0$, or $K - \frac{1}{\omega^2 \rho} \mathrm{id} : U' \to U'$ is a homeomorphism. In the second case this implies $-\omega^2 \rho (K - \frac{1}{\omega^2 \rho} \mathrm{id}) B : U \to U'$ is a homeomorphism, so there exists a unique solution $u \in U$ to the equation $-\omega^2 \rho (K - \frac{1}{\omega^2 \rho} \mathrm{id}) Bu = -\omega^2 \rho \iota u + Bu = \ell$ for any $\ell \in U'$ which satisfies that $\|u\|_U \leq \|(-\omega^2 \rho (K - \frac{1}{\omega^2 \rho} \mathrm{id}) B)^{-1}\| \|\ell\|_{U'}$. When $\omega = 0$ and for any $\ell \in U'$, obviously $\|u\|_U \leq \|B^{-1}\| \|\ell\|_{U'}$, where $u = B^{-1} \ell$ is the unique solution to $Bu = \ell$.

From the theory of compact operators the set of eigenvalues of $K$ is countable, bounded, and can only accumulate at 0. Since the eigenvalues considered are of the form $\frac{1}{\omega^2 \rho}$, it follows that their inverses, $\omega^2 \rho$, are also countable and have no accumulation point.

It remains to show the embedding $\iota : U \to U'$ is compact. This is due to the fact that $(U, V, U')$ is a Gelfand triple, with $V = L^2(\Omega)$. More precisely, the natural embedding $\iota_V : U \to V$, $\iota_V u = u$, is continuous by the Sobolev embedding theorem, and moreover $\overline{U}^V = V$ since $U$ contains all smooth functions vanishing in $\partial \Omega$ which are well known to be dense in $V$. Thus, the transpose $\iota_V^\mathsf{T} : V' \to U'$ is continuous, takes the form $\iota_V^\mathsf{T} v = v|_U$, and is injective by the density of $U$ in $V$. Let $R_V : V \to V'$ be the Riesz map, which explicitly takes the form $\langle R_V u, v \rangle = (u,v)_\Omega$, and is known to be continuous and bijective by the Riesz representation theorem. Thus, the original embedding $\iota_V^\mathsf{T} R_V \iota_V = \iota : U \to U'$ is injective and compact, because $\iota_V$ is compact by the Rellich-Kondrachov theorem. □

**Remark 2.3.** The theorem can be generalized to spatially heterogeneous (but constant in time) densities, as long as $\rho_{\min} < \rho(x) < \rho_{\max}$ for all $x \in \Omega$, where $\rho_{\min} > 0$ and $\rho_{\max} > 0$ are constants.

**Remark 2.4.** Theorem 2.1 shows that (2.5) is well-posed for almost every value of $\omega$, with the exception of some critical values which are essentially spread out in the real-number line. At these critical values the system is said to be in resonance, and a unique solution does not exist. Indeed, the constant $C$ in the statement of the theorem, which is $\omega$-dependent, blows up as these resonant frequencies are approached. Thus, when close to these frequencies, numerical schemes discretizing these equations, even if theoretically stable, are usually very ill-conditioned and round-off error may play an undesirable role (see [45]).

## 2.3 Broken primal formulation

In the study of discontinuous finite element methods it is common to merely consider functions that elementwise have a particular regularity and are possibly discontinuous at the boundaries of the elements, instead of requiring those functions to have the regularity at a global level. This leads to broken spaces dependent on a mesh (a relatively open partition of $\Omega$), $\mathcal{T}$. The broken test space analogous to $H^1(\Omega)$ is

$$H^1(\mathcal{T}) = \{u : \Omega \to \mathbb{C}^3 \mid \|u\|_{H^1(\mathcal{T})}^2 = \sum_{K \in \mathcal{T}} \|u|_K\|_{H^1(K)}^2 < \infty\}. \tag{2.6}$$

Proceeding as with the classical case, but this time multiplying (2.2) by a broken test function $v \in H^1(\mathcal{T})$ yields,

$$-\omega^2 \rho(u,v)_\mathcal{T} + (\mathsf{C}^*{:}\nabla u, \nabla v)_\mathcal{T} - \langle(\mathsf{C}^*{:}\nabla u)\cdot\hat{n}, v\rangle_{\partial\mathcal{T}} = (f,v)_\mathcal{T},$$
$$(u,v)_\mathcal{T} = \sum_{K \in \mathcal{T}} (u|_K, v|_K)_K, \qquad \langle u,v\rangle_{\partial\mathcal{T}} = \sum_{K \in \mathcal{T}} \langle u_K, v_K\rangle_{\partial K}, \quad (2.7)$$

where $\langle\,\cdot\,,\,\cdot\,\rangle_{\partial K}$ for now can be interpreted as a boundary integral in $\partial K$. Unfortunately, in its current form, the formulation in (2.7) is not useful from a theoretical perspective. This issue led to the study of so-called broken variational formulations, which occur in the context of the DPG methodology [6]. The rigorous comprehension of $\langle\,\cdot\,,\,\cdot\,\rangle_{\partial\mathcal{T}}$ is fundamental in developing the broken primal formulation, as well as the introduction of a new interface variable to replace $(\mathsf{C}^*{:}\nabla u)\cdot\hat{n}$ at the skeleton of the mesh (the boundaries of the elements of the mesh), $\partial\mathcal{T}$. Thus, this interface variable is a traction, so intuitively, this amounts to finding the correct functional space for tractions at $\partial\mathcal{T}$, which in turn come from stress fields in $\Omega$.

The appropriate space for stresses is

$$H_{\Gamma_\sigma}(\mathrm{div}, \Omega) = \{\sigma \in H(\mathrm{div}, \Omega) \mid \sigma|_{\Gamma_\sigma}\cdot\hat{n} = 0\}, \quad (2.8)$$

where for any $K \subseteq \Omega$,

$$H(\mathrm{div}, K) = \{\sigma : K \to \mathbb{C}^{3\times 3} \mid \|\sigma\|^2_{H(\mathrm{div},K)} = (\sigma,\sigma)_K + (\mathrm{div}\,\sigma, \mathrm{div}\,\sigma)_K < \infty\},$$
$$H(\mathrm{div}, \mathcal{T}) = \{\sigma : \Omega \to \mathbb{C}^{3\times 3} \mid \|\sigma\|^2_{H(\mathrm{div},\mathcal{T})} = \sum_{K \in \mathcal{T}} \|\sigma|_K\|^2_{H(\mathrm{div},K)} < \infty\}, \quad (2.9)$$

with $\mathrm{div}\,\sigma$ being the row-wise distributional divergence of $\sigma$. Then, the *interface* space for tractions is

$$H^{-1/2}_{\Gamma_\sigma}(\partial\mathcal{T}) = \mathrm{tr}_{2,\mathcal{T}}(H_{\Gamma_\sigma}(\mathrm{div}, \Omega)), \quad (2.10)$$

where

$$\begin{aligned}
\mathrm{tr}_{0,\mathcal{T}} : H^1(\mathcal{T}) \longrightarrow H^{1/2}_\Pi(\partial\mathcal{T}) &= \prod_{K \in \mathcal{T}} H^{1/2}(\partial K) = \prod_{K \in \mathcal{T}} \{\hat{u}_K = u|_{\partial K} \mid u \in H^1(K)\}, \\
\mathrm{tr}_{2,\mathcal{T}} : H(\mathrm{div}, \mathcal{T}) \to H^{-1/2}_\Pi(\partial\mathcal{T}) &= \prod_{K \in \mathcal{T}} H^{-1/2}(\partial K) = \prod_{K \in \mathcal{T}} \{\hat{\sigma}_{n,K} = \sigma|_{\partial K}\cdot\hat{n}_K \mid \sigma \in H(\mathrm{div}, K)\}, \\
\mathrm{tr}_{0,\mathcal{T}} u = \prod_{K \in \mathcal{T}} u|_K\big|_{\partial K}, & \qquad \mathrm{tr}_{2,\mathcal{T}}\sigma = \prod_{K \in \mathcal{T}} \sigma|_K\big|_{\partial K}\cdot\hat{n}_K,
\end{aligned} \quad (2.11)$$

with $\hat{n}_K$ being the outward normal to $K \in \mathcal{T}$. In fact, $H^{1/2}_\Pi(\partial\mathcal{T})$ and $H^{-1/2}_\Pi(\partial\mathcal{T})$ are dual spaces to each other, and so are $H^{1/2}(\partial K)$ and $H^{-1/2}(\partial K)$ for each $K \in \mathcal{T}$ [48]. This gives a rigorous interpretation of $\langle\,\cdot\,,\,\cdot\,\rangle_{\partial\mathcal{T}}$ and $\langle\,\cdot\,,\,\cdot\,\rangle_{\partial K}$, which are duality pairings that become boundary integrals for smooth enough inputs. Finally, it can be shown that $\|\hat{\sigma}_n\|_{H^{-1/2}_\Pi(\partial\mathcal{T})} = \inf_{\sigma \in \mathrm{tr}^{-1}_{2,\mathcal{T}}\{\hat{\sigma}_n\}} \|\sigma\|_{H(\mathrm{div},\mathcal{T})}$ (see [27]), where $\|\cdot\|_{H(\mathrm{div},\mathcal{T})}$ can be replaced by $\|\cdot\|_{H(\mathrm{div},\Omega)}$ when $\hat{\sigma}_n \in H^{-1/2}_{\Gamma_\sigma}(\partial\mathcal{T})$, since $H_{\Gamma_\sigma}(\mathrm{div},\Omega) = \mathrm{tr}^{-1}_{2,\mathcal{T}}(H^{-1/2}_{\Gamma_\sigma}(\partial\mathcal{T}))$.

Assuming vanishing boundary conditions, $u^{\Gamma_u} = 0$ and $\sigma_n^{\Gamma_\sigma} = 0$, the sesquilinear and conjugate-linear forms of the broken variational formulation are

$$b_\mathcal{T}\big((u, \hat{\sigma}_n), v\big) = b_0(u,v) + \hat{b}(\hat{\sigma}_n, v), \qquad \ell_\mathcal{T}(v) = (f,v)_\mathcal{T},$$
$$b_0(u,v) = -\omega^2 \rho(u,v)_\mathcal{T} + (\mathsf{C}^*{:}\nabla u, \nabla v)_\mathcal{T}, \qquad \hat{b}(\hat{\sigma}_n, v) = -\langle\hat{\sigma}_n, \mathrm{tr}_{0,\mathcal{T}}v\rangle_{\partial\mathcal{T}}, \quad (2.12)$$

where $u \in U = H^1_{\Gamma_u}(\Omega)$, $\hat{\sigma}_n \in \hat{U} = H^{-1/2}_{\Gamma_\sigma}(\partial\mathcal{T})$, $v \in V_\mathcal{T} = H^1(\mathcal{T})$, and the trial space $U_\mathcal{T} = U \times \hat{U}$ is equipped with its Hilbert norm. Note that in relation to (2.3), the domain of the test space of $b_0$ was extended from $H^1_{\Gamma_u}(\Omega)$ to $H^1(\mathcal{T})$. Meanwhile, it is clear $\ell_\mathcal{T} \in V'_\mathcal{T}$. Thus, the broken primal formulation that solves (2.2) is,

$$\begin{cases} \text{Find } (u, \hat{\sigma}_n) \in U_\mathcal{T}\,, \\ b_\mathcal{T}\big((u, \hat{\sigma}_n), v\big) = \ell_\mathcal{T}(v)\,, \quad \text{for all } v \in V_\mathcal{T}. \end{cases} \quad (2.13)$$

When the boundary conditions are nontrivial, terms involving extensions of $u^{\Gamma_u} \in \text{tr}_{0,\{\Omega\}}(H^1(\Omega))|_{\Gamma_u}$ and $\sigma_n^{\Gamma_\sigma} \in \text{tr}_{2,\{\Omega\}}(H(\text{div},\Omega))|_{\Gamma_\sigma}$ to $H^1(\Omega)$ and $\text{tr}_{2,\mathcal{T}}(H(\text{div},\Omega))$ respectively, become part of $\ell_\mathcal{T}$.

In [43] it was proved that $\langle \hat{\sigma}_n, \text{tr}_{0,\mathcal{T}} v \rangle_{\partial \mathcal{T}} = 0$ for all $\hat{\sigma}_n \in H^{-1/2}_{\Gamma_\sigma}(\partial \mathcal{T})$ if and only if $v \in V_0 = H^1_{\Gamma_u}(\Omega) \subseteq V_\mathcal{T}$. To begin with, this implies that $\hat{b}|_{\hat{U} \times V_0} = 0$, and that $b_\mathcal{T}|_{U \times V_0}$ and $\ell_\mathcal{T}|_{V_0}$ are effectively the sesquilinear and conjugate-linear forms of the classical primal formulation (since $\hat{U}$ ceases to play a role). Moreover, this fact and [6, Theorem 2.3] yield the well-posedness of the broken primal formulation via a straightforward application of [6, Theorem 3.1] (see also [43, Theorem 3.1]), provided the classical primal formulation is well-posed. Thus, assuming $|\bar{\varepsilon} : \mathfrak{Re}(\mathsf{C}^*) : \varepsilon| > 0$ for all symmetric second-order tensors $\varepsilon \neq 0$, the broken primal formulation is well-posed for most values of $\omega$ as established by Theorem 2.1.

**Theorem 2.2.** *Let $U = H^1_{\Gamma_u}(\Omega)$, $\hat{U} = H^{-1/2}_{\Gamma_\sigma}(\partial \mathcal{T})$, $U_\mathcal{T} = U \times \hat{U}$, $V_\mathcal{T} = H^1(\mathcal{T})$, and consider the problem in (2.13), with $b_\mathcal{T}$ defined in (2.12) in terms of $b_0$ and $\hat{b}$. Then, (2.13) is well-posed if and only if the problem in (2.5) is well-posed. In case of being well-posed, given any $\ell_\mathcal{T} \in V'_\mathcal{T}$, there exists a unique solution $(u, \hat{\sigma}_n) \in U_\mathcal{T}$ solving $b_\mathcal{T}\big((u, \hat{\sigma}_n), v\big) = \ell_\mathcal{T}(v)$ for all $v \in V_\mathcal{T}$ which satisfies $\|(u, \hat{\sigma}_n)\|_{U_\mathcal{T}} \leq C \|\ell_\mathcal{T}\|_{V'_\mathcal{T}}$ for a $C > 0$ independent of the choice of $\ell_\mathcal{T}$ and mesh $\mathcal{T}$.*

## 3 Numerical method

### 3.1 Minimum residual methods

Minimum residual finite element methods begin with the most general setup of an arbitrary linear variational formulation,

$$\begin{cases} \text{Find } \mathfrak{u} \in U\,, \\ b(\mathfrak{u}, \mathfrak{v}) = \ell(\mathfrak{v})\,, \quad \text{for all } \mathfrak{v} \in V\,, \end{cases} \Leftrightarrow \begin{cases} \text{Find } \mathfrak{u} \in U\,, \\ B\mathfrak{u} = \ell\,, \end{cases} \quad (3.1)$$

where $b$ is a sesquilinear form with a Hilbert trial space $U$ and Hilbert test space $V$, $\ell$ is a conjugate-linear form on $V$, and $B : U \to V'$ is a linear operator defined by $\langle B\mathfrak{u}, \mathfrak{v} \rangle = b(\mathfrak{u}, \mathfrak{v}) = \langle B^\dagger \mathfrak{v}, \mathfrak{u} \rangle$, with $B^\dagger : V \to U'$ being the conjugate-linear transpose of $B$. Notice $U'$ is the continuous dual of $U$, but $V'$ is the continuous conjugate-dual of $V$ consisting of conjugate-linear functionals. The next step is to initiate the discretization of the problem by considering a finite-dimensional discrete trial space $U_h \subseteq U$. Then, simply attempt to find the miminizer of the residual $B\mathfrak{u} - \ell$ over $U_h$,

$$\mathfrak{u}_h = \arg\min_{\mathfrak{u} \in U_h} \|B\mathfrak{u} - \ell\|^2_{V'}\,. \quad (3.2)$$

Computing the Gâteaux derivative yields that the solution $\mathfrak{u}_h \in U_h$ satisfies that $(B\mathfrak{u}_h - \ell, B\delta\mathfrak{u})_{V'} = 0$ for all $\delta\mathfrak{u} \in U_h$, which can be recast as

$$\begin{cases} \text{Find } \mathfrak{u}_h \in U_h\,, \\ b(\mathfrak{u}_h, \delta\mathfrak{v}) = \ell(\delta\mathfrak{v})\,, \quad \text{for all } \delta\mathfrak{v} = R_V^{-1} B\delta\mathfrak{u} \in R_V^{-1} BU_h = V^{\text{opt}}\,, \end{cases} \Leftrightarrow \begin{cases} \text{Find } \mathfrak{u}_h \in U_h\,, \\ B^\dagger R_V^{-1} B\mathfrak{u}_h = B^\dagger R_V^{-1} \ell\,, \end{cases} \quad (3.3)$$

where $R_V : V \to V'$ is the Riesz map of $V$, which is defined by $\langle R_V \mathfrak{v}, \delta\mathfrak{v}\rangle_{V'\times V} = (\mathfrak{v}, \delta\mathfrak{v})_V$ for all $\mathfrak{v}, \delta\mathfrak{v} \in V$. Here, $V^{\mathrm{opt}} = R_V^{-1} B U_h$ is called the optimal test space, and clearly satisfies that $\dim(U_h) = \dim(V^{\mathrm{opt}})$, since the Riesz map is an isometric isomorphism. Moreover, by isometry, $\|B\mathfrak{u}\|_{V'} = \|R_V^{-1} B\mathfrak{u}\|_V$, and the optimality becomes clear as it can be shown that

$$\gamma^{\mathrm{opt}} = \inf_{\mathfrak{u}_h \in U_h} \sup_{\delta\mathfrak{v} \in V^{\mathrm{opt}}} \frac{|b(\mathfrak{u}_h, \delta\mathfrak{v})|}{\|\mathfrak{u}_h\|_U \|\delta\mathfrak{v}\|_V} = \inf_{\mathfrak{u}_h \in U_h} \sup_{\mathfrak{v} \in V} \frac{|b(\mathfrak{u}_h, \mathfrak{v})|}{\|\mathfrak{u}_h\|_U \|\mathfrak{v}\|_V} \geq \inf_{\mathfrak{u} \in U} \sup_{\mathfrak{v} \in V} \frac{|b(\mathfrak{u}, \mathfrak{v})|}{\|\mathfrak{u}\|_U \|\mathfrak{v}\|_V} = \gamma, \qquad (3.4)$$

where the infima and suprema are tacitly assumed to be taken over nonzero elements. Naturally, the original problem in (3.1) is assumed to be well-posed, implying $\gamma^{\mathrm{opt}} \geq \gamma > 0$, and by Babuška's theorem [2], the problem in (3.3) is said to be stable, so that there exists a unique solution $\mathfrak{u}_h \in U_h$ that satisfies the stability estimate, $\|\mathfrak{u}_h\|_U \leq \frac{1}{\gamma}\|\ell\|_{(V^{\mathrm{opt}})'}$, as well as the relation [2, 68] (see [64, 49] for the Banach space setting),

$$\|\mathfrak{u} - \mathfrak{u}_h\|_U \leq \frac{\|b\|}{\gamma} \inf_{\delta\mathfrak{u}_h \in U_h} \|\mathfrak{u} - \delta\mathfrak{u}_h\|_U, \qquad (3.5)$$

where $\|b\| = \sup_{(\mathfrak{u},\mathfrak{v}) \in U \times V} \frac{|b(\mathfrak{u},\mathfrak{v})|}{\|\mathfrak{u}\|_U \|\mathfrak{v}\|_V}$, and $\mathfrak{u} \in U$ is the unique solution to (3.1). This setting is referred to as the ideal or optimal Petrov-Galerkin method.

Regrettably, this ideal method is not computationally viable in most cases, since the Riesz map $R_V$ cannot be inverted as $V$ is an infinite-dimensional space. Instead, minimum residual methods invert the Riesz map over a finite-dimensional enriched test space $V^{\mathrm{enr}} \subseteq V$ satisfying that $\dim(V^{\mathrm{enr}}) \geq \dim(U_h)$. Therefore, (3.2) is minimized with the norm $\|\cdot\|_{(V^{\mathrm{enr}})'}$ instead, and equivalently $R_V$ is replaced with $R_{V^{\mathrm{enr}}}$ in (3.3), so that the new test space is $V_h = R_{V^{\mathrm{enr}}}^{-1} B U_h$. Clearly, $V_h$ aims to approximate the optimal test space $V^{\mathrm{opt}}$, and intuitively, the larger $V^{\mathrm{enr}}$ is, the closer $V_h$ will be to $V^{\mathrm{opt}}$. From a computational standpoint it is more convenient to consider the second characterization in (3.3) which is in $U'$ and leads to

$$B^\dagger R_{V^{\mathrm{enr}}}^{-1} B \mathfrak{u}_h = B^\dagger R_{V^{\mathrm{enr}}}^{-1} \ell. \qquad (3.6)$$

From the point of view of linear algebra, this equation can be discretized as

$$\mathsf{B}^{\mathrm{n\text{-}opt}} \mathsf{u}_h = \mathsf{B}^\dagger \mathsf{G}^{-1} \mathsf{B} \mathsf{u}_h = \mathsf{B}^\dagger \mathsf{G}^{-1} \mathsf{l} = \mathsf{l}^{\mathrm{n\text{-}opt}}, \qquad (3.7)$$

where $\mathsf{B}_{ij} = b(\mathfrak{u}_j, \mathfrak{v}_i)$, $\mathsf{l}_i = \ell(\mathfrak{v}_i)$, $\mathsf{G}_{ij} = (\mathfrak{v}_i, \mathfrak{v}_j)_V$, and $\mathfrak{u}_h = \sum_{j=1}^{\dim(U_h)} (\mathsf{u}_h)_j \mathfrak{u}_j$, with $\{\mathfrak{u}_j\}_{j=1}^{\dim(U_h)}$ and $\{\mathfrak{v}_i\}_{i=1}^{\dim(V^{\mathrm{enr}})}$ being bases for $U_h$ and $V^{\mathrm{enr}}$ respectively. Here, $\mathsf{B}^\dagger$ is the Hermitian transpose of the tall rectangular matrix $\mathsf{B}$, $\mathsf{G}$ is a Gram matrix, and $\mathsf{B}^{\mathrm{n\text{-}opt}} = \mathsf{B}^\dagger \mathsf{G}^{-1} \mathsf{B}$ and $\mathsf{l}^{\mathrm{n\text{-}opt}} = \mathsf{B}^\dagger \mathsf{G}^{-1} \mathsf{l}$ are called the near-optimal stiffness matrix and load. First, notice that via this procedure,

$$V_h = \mathrm{span}\left(\left\{\sum_{i=1}^{\dim(V^{\mathrm{enr}})} (\mathsf{G}^{-1}\mathsf{B})_{ij} \mathfrak{v}_i\right\}_{j=1}^{\dim(U_h)}\right), \qquad (3.8)$$

so $V_h$ is inherently being computed from $V^{\mathrm{enr}}$, with the big advantage that the basis $\{\mathfrak{v}_i\}_{i=1}^{\dim(V^{\mathrm{enr}})}$ may be a standard discretization of $V$. This distinguishes the method from other Petrov-Galerkin methods, where finding an exotic basis for $V_h$ is typically required from the start, but here only a standard basis of $V^{\mathrm{enr}}$ (not $V_h$) is sufficient. Second, notice the stiffness matrix, $\mathsf{B}^{\mathrm{n\text{-}opt}} = (\mathsf{G}^{-1/2}\mathsf{B})^\dagger (\mathsf{G}^{-1/2}\mathsf{B})$, is always Hermitian and positive definite.

Finally, it is worth pointing out that the error estimate in (3.5) will no longer hold since $V_h \neq V^{\mathrm{opt}}$. However, using a Fortin operator it can be shown that a similar estimate still holds. Assume the existence of a Fortin operator, $\Pi_F : V \to V^{\text{enr}}$, defined such that it is linear, continuous and satisfies the orthogonality condition $b(\mathfrak{u}_h, \mathfrak{v} - \Pi_F \mathfrak{v}) = 0$ for all $\mathfrak{u}_h \in U_h$ and $\mathfrak{v} \in V$ [32, 50]. Then, it can be shown that the error estimate becomes

$$\|\mathfrak{u} - \mathfrak{u}_h\|_U \leq \frac{\|b\| M_F}{\gamma} \inf_{\delta \mathfrak{u}_h \in U_h} \|\mathfrak{u} - \delta \mathfrak{u}_h\|_U, \qquad (3.9)$$

where $M_F \geq \|\Pi_F\| = \sup_{\mathfrak{v} \in V} \frac{\|\Pi_F \mathfrak{v}\|_V}{\|\mathfrak{v}\|_V}$. Similarly, the stability estimate of the solution $\mathfrak{u}_h \in U_h$ to (3.7) becomes $\|\mathfrak{u}_h\|_U \leq \frac{M_F}{\gamma} \|\ell\|_{V_h'}$. Note that Fortin operators only yield conservative estimates of the actual values of the constants.

## 3.2 The DPG methodology

Unfortunately, the problem in (3.7) may still be computationally prohibitive as it requires computing $\mathsf{G}^{-1}$ beforehand, and this is a global problem which may be expensive. The solution is to consider variational formulations with broken test spaces, which are referred to as broken variational formulations. The application of minimum residual methods to broken variational formulations is called the DPG methodology. This allows to localize computations, including the inversion of the Riesz map which can be done elementwise. Moreover, the construction of Fortin operators can be made local too, and indeed these can typically be constructed for large enough $V^{\text{enr}}$ as in [32, 6, 50].

More specifically, for broken variational formulations the discrete trial space and enriched test space take the form,

$$\begin{aligned} U_h &= \{\delta \mathfrak{u}_h \in U \mid \delta \mathfrak{u}_{h,K} \in U_h(K), \forall K \in \mathcal{T}\}, \\ V^{\text{enr}} &= \{\mathfrak{v}^{\text{enr}} \mid \mathfrak{v}^{\text{enr}}|_K \in V^{\text{enr}}(K), \forall K \in \mathcal{T}\} \cong \bigoplus_{K \in \mathcal{T}} V^{\text{enr}}(K), \end{aligned} \qquad (3.10)$$

where $\mathfrak{v}^{\text{enr}}|_K$ represents the restriction of the domain from $\Omega$ to $K \in \mathcal{T}$, and where $\delta \mathfrak{u}_{h,K}$, $U_h(K)$, $V^{\text{enr}}(K)$ will be defined later. The last congruence follows because $V^{\text{enr}}$ is locally decoupled, and in fact its basis may be written as

$$\{\mathfrak{v}_i\}_{i=1}^{\dim(V^{\text{enr}})} = \bigcup_{K \in \mathcal{T}} \{\mathfrak{v}_{i_K}\}_{i_K=1}^{\dim(V^{\text{enr}}(K))}, \qquad (3.11)$$

so that for each $\mathfrak{v}_j$ there is a *unique* $K \in \mathcal{T}$ for which $\mathfrak{v}_j|_K \neq 0$. Since $\mathsf{G}_{ij} = (\mathfrak{v}_i, \mathfrak{v}_j)_V$, it follows $\mathsf{G}$ will have a block-diagonal structure when its basis is organized and indexed by $K \in \mathcal{T}$, where $\mathsf{G}_K$ is the $K$-th block. Thus, $\mathsf{G}^{-1}$ will also be block-diagonal with its $K$-th block being $\mathsf{G}_K^{-1}$. In fact, the size of $\mathsf{G}_K$ is reasonable enough such that a direct Cholesky factorization, $\mathsf{G}_K = \mathsf{L}_K^\dagger \mathsf{L}_K$, is viably computed. Hence, the stiffness matrix and load in (3.7) can be assembled from its local versions as usual,

$$\begin{aligned} \mathsf{B}^{\text{n-opt}} &= \sum_{K \in \mathcal{T}} \mathsf{A}_K^\mathsf{B} \mathsf{B}_K^{\text{n-opt}}, & \mathsf{B}_K^{\text{n-opt}} &= \mathsf{B}_K^\dagger \mathsf{G}_K^{-1} \mathsf{B}_K = (\mathsf{L}_K^{-1} \mathsf{B}_K)^\dagger (\mathsf{L}_K^{-1} \mathsf{B}_K), \\ \mathsf{l}^{\text{n-opt}} &= \sum_{K \in \mathcal{T}} \mathsf{A}_K^\mathsf{l} \mathsf{l}_K^{\text{n-opt}}, & \mathsf{l}_K^{\text{n-opt}} &= \mathsf{B}_K^\dagger \mathsf{G}_K^{-1} \mathsf{l}_K = (\mathsf{L}_K^{-1} \mathsf{B}_K)^\dagger \mathsf{L}_K^{-1} \mathsf{l}_K, \end{aligned} \qquad (3.12)$$

where $\mathsf{A}_K^\mathsf{B}$ and $\mathsf{A}_K^\mathsf{l}$ are standard local to global assembly operators.

Next, a natural a posteriori error estimator for these methods will be described. It is said to be natural, because most a posteriori estimators are based on the residual, and here, the methods themselves are designed to minimize the residual as expressed in (3.2). Indeed, the expression for the minimum residual in (3.2) is given by $\|B\mathfrak{u}_h - \ell\|_{V'}^2 = \langle B\mathfrak{u}_h - \ell, R_V^{-1}(B\mathfrak{u}_h - \ell) \rangle$. In practice, only $V^{\text{enr}}$ can be used for $V$, so the residual is approximated by $\|B\mathfrak{u}_h - \ell\|^2_{(V^{\text{enr}})'} = (B\mathfrak{u}_h - \mathfrak{l})^\dagger \mathsf{G}^{-1}(B\mathfrak{u}_h - \mathfrak{l})$ instead. With broken variational formulations, this expression can be localized to an element residual for each $K \in \mathcal{T}$,

$$\mathsf{r}_K^2 = (\mathsf{B}_K \mathsf{u}_{h,K} - \mathsf{l}_K)^\dagger \mathsf{G}_K^{-1}(\mathsf{B}_K \mathsf{u}_{h,K} - \mathsf{l}_K) = (\mathsf{L}_K^{-1}(\mathsf{B}_K \mathsf{u}_{h,K} - \mathsf{l}_K))^\dagger (\mathsf{L}_K^{-1}(\mathsf{B}_K \mathsf{u}_{h,K} - \mathsf{l}_K)). \tag{3.13}$$

This serves as a natural a posteriori error estimator to drive adaptivity and is valid for any arbitrary polynomial order $p$. Indeed, by construction, minimum residual methods attempt to reduce the residual $\|B\mathfrak{u}_h - \ell\|_{(V^{\text{enr}})'}$ over $U_h$, so this value is always expected to go down globally as the mesh is refined (which is guaranteed under $p$ refinements of $U_h$, but not necessarily under $h$ refinements as the interface spaces do not embed when the mesh is refined).

## 3.3 Convergence analysis and exact sequence spaces

**General arguments for $h$-convergence**

In this section, the convergence of the numerical method will be analyzed for linear viscoelasticity. For now assume the existence of a Fortin operator so that the stability error estimate in (3.9) holds. For $U_\mathcal{T} = H^1_{\Gamma_u}(\Omega) \times H^{-1/2}_{\Gamma_\sigma}(\partial \mathcal{T})$ in (2.13), the bound explicitly takes the form,

$$\|u - u_h\|^2_{H^1(\Omega)} + \|\hat{\sigma}_n - \hat{\sigma}_{n,h}\|^2_{H^{-1/2}_\Pi(\partial \mathcal{T})} \leq C_{\text{st}}^2 \inf_{(w_h, \hat{\tau}_{n,h}) \in U_h} \left( \|u - w_h\|^2_{H^1(\Omega)} + \|\hat{\sigma}_n - \hat{\tau}_{n,h}\|^2_{H^{-1/2}_\Pi(\partial \mathcal{T})} \right), \tag{3.14}$$

where the stability constant is $C_{\text{st}} = \frac{\|b_\mathcal{T}\| M_F}{\gamma}$.

To proceed further, the discrete trial space, $U_h$, and enriched test space, $V_\mathcal{T}^{\text{enr}}$, must be defined. For this, consider an affine shape-regular element $K \in \mathcal{T}$ and an arbitrary polynomial order $p$. Assume there exist conforming discretizations for $H^1(K)$, $H(\text{curl}, K)$, $H(\text{div}, K)$ and $L^2(K)$, satisfying the exact sequence property and containing the appropriate polynomials,

$$\begin{array}{cccc}
H^1(K) & H(\text{curl}, K) & H(\text{div}, K) & L^2(K) \\
\cup\mathsf{I} & \cup\mathsf{I} & \cup\mathsf{I} & \cup\mathsf{I} \\
W^p(K) \xrightarrow{\nabla} & Q^p(K) \xrightarrow{\nabla \times} & V^p(K) \xrightarrow{\nabla \cdot} & Y^p(K) \\
\cup\mathsf{I} & \cup\mathsf{I} & \cup\mathsf{I} & \cup\mathsf{I} \\
\mathcal{P}^p & (\mathcal{P}^{p-1})^3 & (\mathcal{P}^{p-1})^3 & \mathcal{P}^{p-1},
\end{array} \tag{3.15}$$

where $\mathcal{P}^p$ are the set of polynomials in three variables of total order $p$. Then, for any $s > \frac{1}{2}$, the projection-based interpolation operators, $\Pi_i^{K,s}$ for $i = 0, 1, 2, 3$, generate a commuting exact sequence [16],

$$\begin{array}{ccccc}
H^{1+s}(K) & \xrightarrow{\nabla} & H^s(\text{curl}, K) & \xrightarrow{\nabla \times} & H^s(\text{div}, K) & \xrightarrow{\nabla \cdot} & H^s(K) \\
\downarrow \Pi_0^{K,s} & & \downarrow \Pi_1^{K,s} & & \downarrow \Pi_2^{K,s} & & \downarrow \Pi_3^{K,s} \\
W^p(K) & \xrightarrow{\nabla} & Q^p(K) & \xrightarrow{\nabla \times} & V^p(K) & \xrightarrow{\nabla \cdot} & Y^p(K),
\end{array} \tag{3.16}$$

where the fractional $s > \frac{1}{2}$ Sobolev spaces are slightly more regular counterparts of the usual Sobolev spaces in (3.15) (which correspond to $s = 0$) [48]. Here, there is the tacit assumption that the operators are in fact pullbacks of operators associated to scaled master elements of $K$.

The discrete trial space is chosen as

$$\begin{aligned}
U_h &= \left\{ \delta \mathfrak{u}_h = (w_h, \hat{\tau}_{n,h}) \in U_\mathcal{T} \mid \delta \mathfrak{u}_{h,K} = (w_h|_K, \hat{\tau}_{n,h,K}) \in U_h(K), \forall K \in \mathcal{T} \right\}, \\
U_h(K) &= (W^p(K))^3 \times (V^p(\partial K))^3, \qquad (V^p(\partial K))^3 = \{\tau_h|_{\partial K} \cdot \hat{n}_K \mid \tau_h \in (V^p(K))^3\},
\end{aligned} \tag{3.17}$$

where $U_\mathcal{T} = H^1_{\Gamma_u}(\Omega) \times H^{-1/2}_{\Gamma_\sigma}(\partial\mathcal{T})$. Note that the condition $\delta\mathfrak{u}_h = (w_h, \hat{\tau}_{n,h}) \in U_\mathcal{T}$ implies that $w_h$ and $\hat{\tau}_{n,h}$ vanish at $\Gamma_u$ and $\Gamma_\sigma$ respectively, and that $w_h|_{K_1}|_e = w_h|_{K_2}|_e$ and $\hat{\tau}_{n,h,K_1}|_e = -\hat{\tau}_{n,h,K_2}|_e$ where $e$ is an edge shared by $K_1$ and $K_2$. Hence, the definition of $U_h$ implies that one must ensure compatibility across elements. Meanwhile, the enriched test space is selected from a *sequence* of order $p + \Delta p$, so that

$$V^{\text{enr}}_\mathcal{T} = \left\{\mathfrak{v}^{\text{enr}} = v \mid \mathfrak{v}^{\text{enr}}|_K = v|_K \in V^{\text{enr}}_\mathcal{T}(K), \forall K \in \mathcal{T}\right\} \subseteq V_\mathcal{T}, \qquad V^{\text{enr}}_\mathcal{T}(K) = (W^{p+\Delta p}(K))^3, \qquad (3.18)$$

where $V_\mathcal{T} = H^1(\mathcal{T})$.

Next, consider $\hat{\eta}_{n,K} = \eta|_{\partial K} \cdot \hat{n}_K$ for some $\eta \in H^s(\text{div}, K)$ and its minimum energy extension norm, so that $\|\hat{\eta}_{n,K}\|_{H^{-1/2}(\partial K)} \leq \|\eta\|_{H(\text{div},K)}$, and

$$\|\hat{\eta}_{n,K} - (\Pi^{K,s}_2 \eta)|_{\partial K} \cdot \hat{n}_K\|_{H^{-1/2}(\partial K)} = \|(\eta - \Pi^{K,s}_2 \eta)|_{\partial K} \cdot \hat{n}_K\|_{H^{-1/2}(\partial K)} \leq \|\eta - \Pi^{K,s}_2 \eta\|_{H(\text{div},K)}. \qquad (3.19)$$

Given that the spaces are assumed to contain the relevant polynomials (see (3.15)), that the sequences commute in (3.16), that the elements are assumed to be affine shape-regular (among all meshes considered), and using standard scaling arguments [24] along with (3.19), yields the usual $h$-convergence interpolation estimates (for any fixed $p$) in $K$ for any sufficiently regular $u$ and $\hat{\sigma}_{n,K}$. The estimates hold globally in $\Omega$ too, since the interpolants give the necessary compatibility across the elements. Thus, choosing the interpolant in (3.14) yields the $h$-convergence result,

$$\left(\|u - u_h\|^2_{H^1(\Omega)} + \|\hat{\sigma}_n - \hat{\sigma}_{n,h}\|^2_{H^{-1/2}_\Pi(\partial\mathcal{T})}\right)^{1/2} \leq C_{\text{st}} C_s h^{\min\{s,p\}} \left(\|u\|^2_{H^{1+s}(\Omega)} + \|\hat{\sigma}_n\|^2_{H^{-1/2+s}(\partial\mathcal{T})}\right)^{1/2}, \qquad (3.20)$$

where $(u, \hat{\sigma}_n) \in U^s_\mathcal{T} \subseteq U_\mathcal{T}$ is the exact solution of (2.13), $(u_h, \hat{\sigma}_{n,h}) \in U_h$ is the computed solution from the DPG methodology, $C_{\text{st}} = \frac{\|b_\mathcal{T}\| M_F}{\gamma}$ is the stability constant, $C_s$ is the $h$-interpolation constant, and $h = \sup_{K \in \mathcal{T}} \text{diam}(K)$. The constant $C_s$ is dependent on $s$, $p$, and the shape-regularity of the elements, but not on $h$. As usual, the $h$-convergence is dictated by the polynomial order or the regularity of the exact solution, i.e., the largest $s$ for which $(u, \hat{\sigma}_n) \in U^s_\mathcal{T}$, where $U^s_\mathcal{T} = H^{1+s}(\Omega) \times H^{-1/2+s}(\partial\mathcal{T})$ has the norm in (3.20). This result is now summarized in the following theorem.

**Theorem 3.1.** *Let $p \in \mathbb{N}$ and take a set of meshes whose elements satisfy a robust affine shape-regularity condition. Consider the problem of finding $\mathfrak{u}_h \in U_h$, such that*

$$b_\mathcal{T}(\mathfrak{u}_h, \mathfrak{v}_h) = \ell_\mathcal{T}(\mathfrak{v}_h), \qquad \forall \mathfrak{v}_h \in V_h,$$

*where $b_\mathcal{T}$ and $\ell_\mathcal{T}$ are defined in (2.12), $U_h \subseteq U_\mathcal{T}$ is defined in (3.17) (assuming the properties in (3.15) hold), $V^{\text{enr}}_\mathcal{T} \subseteq V_\mathcal{T}$ is defined in (3.18), and $V_h$ is defined in terms of a basis of $V^{\text{enr}}_\mathcal{T}$ in (3.8). If a linear and continuous Fortin operator, $\Pi_F : V_\mathcal{T} \to V^{\text{enr}}_\mathcal{T}$, satisfying the orthogonality condition, $b(\mathfrak{u}_h, \mathfrak{v} - \Pi_F \mathfrak{v}) = 0$, for all $\mathfrak{u}_h \in U_h$ and $\mathfrak{v} \in V_\mathcal{T}$ exists, and if $\omega$ is not in the set of critical set of values making (2.13) ill-posed (see Theorem 2.1), then the problem has a unique solution $\mathfrak{u}_h \in U_h$. The discrete solution may be computed using the linear system in (3.7) and the simplifications in (3.12). The exact solution of the infinite-dimensional problem in (2.13) is $\mathfrak{u} \in U_\mathcal{T}$. Provided a $\Pi_F$ continuity bound independent of the meshes, and that $\mathfrak{u}$ is regular enough ($\mathfrak{u} \in U^s_\mathcal{T} \subseteq U_\mathcal{T}$ for $s > \frac{1}{2}$) the discrete solution relates to the exact solution by,*

$$\|\mathfrak{u} - \mathfrak{u}_h\|_{U_\mathcal{T}} \leq C h^{\min\{s,p\}} \|\mathfrak{u}\|_{U^s_\mathcal{T}},$$

*where $h = \sup_{K \in \mathcal{T}} \text{diam}(K)$ and $C = C(s, p, \Omega) > 0$ is independent of the mesh being considered.*

**Verification of the assumptions**

The assumptions in Theorem 3.1 actually do hold. In fact, arbitrary-$p$ discretizations of the type in (3.15) exist for any affinely transformed hexahedron, tetrahedron, triangular prism, and pyramid, so, forgetting about the Fortin operator, the convergence result in Theorem 3.1 would hold for any affine shape-regular hybrid mesh. These discretizations can be found in the literature [28, 10]. For computations, it is fundamental to find explicit expressions for basis functions of these discretizations which are compatible across the different elements. Such hierarchical bases have been presented in [28].

It remains to establish the existence of a Fortin operator. It can be constructed locally due to the broken test spaces. This will be done only for the tetrahedron and with the discrete spaces in (3.15) coming from the Nédélec sequence of the first type [51] (which is used in [28]). Therefore, the local discrete trial space in (3.17), $U_h(K)$, and the local enriched test space in (3.18), $V_{\mathcal{T}}^{\mathrm{enr}}(K)$, are defined by the affine-invariant spaces

$$W^p(K) = \mathcal{P}^p, \qquad V^p(K) = \mathcal{RT}^p = (\mathcal{P}^{p-1})^3 + x\mathcal{P}^{p-1}, \tag{3.21}$$

where $\mathcal{RT}^p$ is also known as the Raviart-Thomas space. Then, [6, Theorem 5.1] establishes a local Fortin operator $\Pi_F^K = \Pi_{F,0}^{K,p,\Delta p} : V_{\mathcal{T}}(K) \to V_{\mathcal{T}}^{\mathrm{enr}}(K)$, which is linear, bounded and satisfies,

$$(\psi_h, \Pi_F^K v - v)_K = 0, \qquad (\phi_h, \nabla(\Pi_F^K v - v))_K = 0, \qquad \langle \hat{\eta}_{n,h}, (\Pi_F^K v - v)|_{\partial K} \rangle_{\partial K} = 0, \tag{3.22}$$

for all $v \in V_{\mathcal{T}}(K) = H^1(K)$, $\psi_h \in (\mathcal{P}^{p+\Delta p-4})^3$, $\phi_h \in (\mathcal{P}^{p+\Delta p-3})^{3\times 3}$ and $\hat{\eta}_{n,h} \in V^{p+\Delta p-2}(\partial K)$, where $\Delta p \geq 3$ and $p \geq 1$. Let $(w_h, \hat{\tau}_{n,h}) \in U_h(K)$. If $\mathsf{C}^*$ is piecewise constant with respect to the mesh, then $w_h \in (\mathcal{P}^p)^3$, $\mathsf{C}^*{:}\nabla w_h \in (\mathcal{P}^{p-1})^{3\times 3}$, and $\hat{\tau}_{n,h} \in V^p(\partial K)$. Thus, provided $\Delta p \geq 4$, it is clear that all the terms in the bilinear form in (2.12) will vanish. For a full tetrahedral mesh, this will hold globally, implying that $b_{\mathcal{T}}((w_h, \hat{\tau}_{n,h}), \Pi_F v - v) = 0$. The bound of the local Fortin operator, $\|\Pi_F^K\|$, is easily seen to be independent of the element size through scaling arguments, so that there does exist a mesh-independent bound $M_F > 0$ such that $\|\Pi_F v\|_V \leq M_F \|v\|_V$ for all $v \in V_{\mathcal{T}}$, provided $\Delta p \geq 4$. However, in principle, the bound $M_F > 0$ may be dependent on the choice of $p$ and $\Delta p$. Proofs of local high-order Fortin operators for the remaining element types (or which are valid for $\Delta p \geq 1$ instead of $\Delta p \geq 4$) and their respective discretizations have not yet been developed, but numerical results suggest they might exist [43]. The following corollary summarizes the result.

**Corollary 3.1.** *All the hypothesis of Theorem 3.1 are satisfied for shape-regular tetrahedral meshes, provided the spaces $W^p(K)$ and $V^p(K)$ in the definitions of $U_h(K)$ and $V_{\mathcal{T}}^{\mathrm{enr}}(K)$ (see (3.17) and (3.18)) are chosen from (3.21), that $\mathsf{C}^*$ is piecewise constant with respect to all the meshes being considered, and that $\Delta p \geq 4$ in (3.18).*

**Remark 3.1.** Due to the nature of the equations, and more specifically to the dynamic term $(u, v)_{\mathcal{T}}$ in (2.12), the requirement of $\Delta p \geq 4$ is more stringent than that proved for linear elasticity and Poisson's equation, which is $\Delta p \geq 3$ [32].

**Comments on $p$- and $hp$-convergence**

Finally, a brief analysis of $p$-convergence is included here as well. For this, given an element $K \in \mathcal{T}$, the existence of a sequence of commuting and bounded polynomial-preserving extension operators from $\partial K$ to $K$ will be assumed, where the bounds must be independent of the polynomial order $p$ [16]. These do exist for tetrahedra [20, 21, 22] and hexahedra [14] with Nédélec's sequences of the first type (like the ones used in [28]), but constructions for the triangular prism and pyramid are still missing. Then, in the $p$-asymptotic limit and for a fine enough mesh, the theory of projection-based interpolation establishes convergence bounds in terms of $p$ as well [16, Theorem 5.3], so that the interpolation constant in (3.20) becomes,

$$C_s = \widetilde{C}_s (\ln p)^2 p^{-s}, \tag{3.23}$$

where $\widetilde{C}_s$ is independent of $p$ and $h$. To finalize the full $hp$-convergence analysis, the only requirement missing is that the Fortin operator bound $M_F$ sitting inside the stability constant $C_{\mathrm{st}} = \frac{\|b_\mathcal{T}\| M_F}{\gamma}$ must be independent of $p$. This has not been established theoretically, but 2D numerical experiments geared directly at local Fortin operators have shown $p$-independence of $M_F$ provided $\Delta p$ is large enough [50]. In any case, the $p$-convergence bounds in (3.23) and the value of $M_F$ are conservative bounds, and in practice the results are much better than the theory predicts.

**Corollary 3.2.** *In Theorem 3.1, if the elements in the meshes are shape-regular tetrahedra or hexahedra, and if there exists a continuity bound of the Fortin operator (an $M_F > 0$ such that $\|\Pi_F \mathfrak{v}\|_{V_\mathcal{T}} \leq M_F \|\mathfrak{v}\|_{V_\mathcal{T}}$ for all $\mathfrak{v} \in V_\mathcal{T}$) which is additionally independent of $p$, then, asymptotically, there exists convergence behavior of the type,*

$$\|\mathfrak{u} - \mathfrak{u}_h\|_{U_\mathcal{T}} \leq \widetilde{C}_s (\ln p)^2 \frac{h^{\min\{s,p\}}}{p^s} \|\mathfrak{u}\|_{U_\mathcal{T}^s},$$

*where $\widetilde{C}_s = \widetilde{C}_s(s, \Omega)$ is independent of both $p$ and $h$.*

# 4 Results

There are a few software libraries that are able to solve DPG methods, like Camellia [57, 58] and DUNE-DPG [33]. In this work, the numerical method described was implemented using an in-house code called *hp3d*, which has support for both local $h$ and $p$ refinements in 3D [24], uses the exact sequence shape functions in [28], and utilizes projection-based interpolation to appropriately enforce the boundary conditions [16]. Thus, the hypotheses of much of the theoretical results in Section 3.3 are satisfied. In the results shown here, only hexahedral and tetrahedral elements were used. The direct solver MUMPS 5.0.1 with OpenMP was used on the global linear system.

First, verification studies confirming the convergence theory were done in a cube. Then, a validation study was completed using results from dynamic mechanical analysis (DMA) calibration experiments on different viscoelastic polymers.

## 4.1 Code verification

To verify the convergence results, a cube, $\Omega = (0, 1)^3$, was discretized initially with five tetrahedra. A manufactured smooth solution for the displacement, $u_i(x) = \prod_{k=1}^{3} \sin(\pi x_k)$, $\forall i = 1, 2, 3$, was utilized to determine the stress, force and boundary data, where the dynamic stiffness tensor, $\mathsf{C}^*$, was defined by $\lambda^* = \mu^* = 1 + \mathrm{i}$. The results are shown in Figure 4.1, where $\Delta p = 1$ in all cases.

Clearly, uniform mesh refinements confirm the $h$-convergence theoretical result in (3.20), since the rate of convergence is of the type $h^p$ due to the analyticity of the solution (so $s = \infty$). When $p \geq s$, where $s$ is the regularity of the solution, Corollary 3.2 establishes an asymptotic $hp$-convergence estimate of the form $\|\mathfrak{u} - \mathfrak{u}_h\|_U \leq C(\ln p)^2 (\frac{h}{p})^s$, where $C$ is independent of $h$ and $p$. This is a quasi-algebraic form of convergence. However, when the solution is analytic, this estimate is expected to improve in some sense, but the explicit form cannot be deduced from (3.23) since $\widetilde{C}_s$ depends on $s$ and its behavior is unknown as $s \to \infty$. Figure 4.1 seems to indicate an exponential $p$-convergence of the form $\|\mathfrak{u} - \mathfrak{u}_h\|_U \leq C \exp(-b\,p^{1.25})$, where $C$ and $b$ are independent of $p$ (but not of $h$) and where $b > 0$ is larger if the mesh is finer. This result can be compared with exponential convergence results found in the literature [39, 40, 62] (it is also better than related $hp$-exponential rates in geometric meshes [3, 34, 35, 61, 60]).

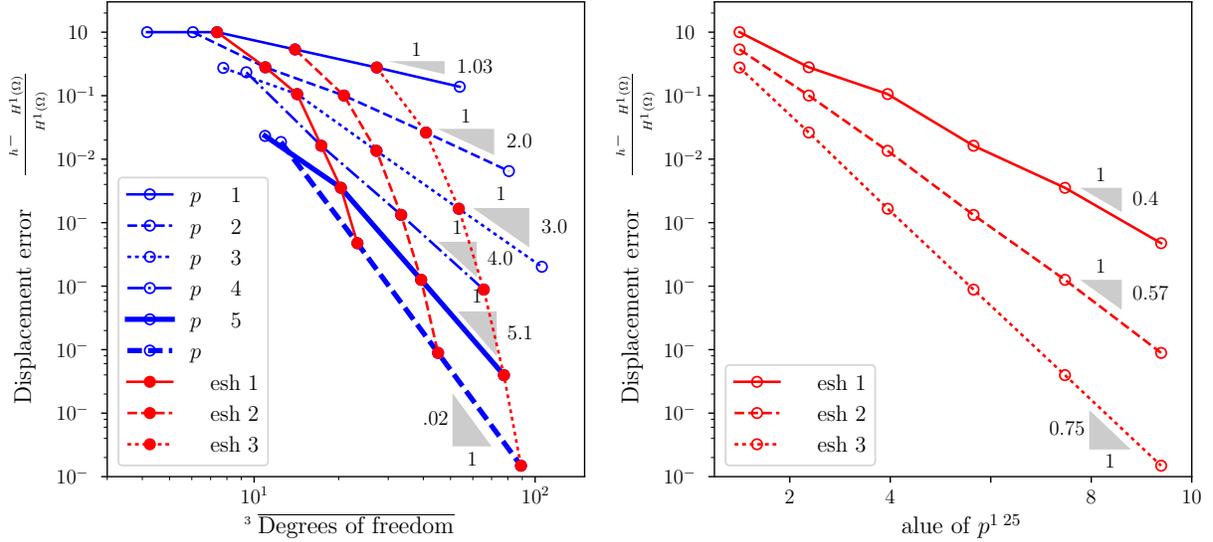

Figure 4.1: Relative displacement error in the $H^1(\Omega)$ norm. Uniform $h$-refinements yield the expected $h^p$ convergence rates for $1 \leq p \leq 6$. Moreover, $p$-refinements of the same mesh do exhibit exponential convergence of the form $\exp(-b\,p^{1.25})$ with $b > 0$ depending on the mesh (the finer the mesh, the higher the $b$).

It should be noted that $\Delta p = 1$ was used in the computations, but numerical experiments were done with higher values of $\Delta p$ as well (including $\Delta p = 4$), and the resulting data points were almost exactly the same. Thus, for this particular equation it seems $\Delta p = 1$ is preferable, since the results are the same and the local computational cost is much lower. However, this merits further theoretical study to be certain, perhaps by finding a Fortin operator which is valid for $\Delta p \geq 1$. Having said that, there are equations and solution schemes where higher values of $\Delta p$ provide advantages (see [15]), so this possibility should not be discarded either.

## 4.2 Validation of DMA experiments

Characterization of viscoelastic material properties in the frequency domain is done through dynamic mechanical analysis (DMA) experiments, where the material is subjected to oscillations. More precisely, to find the dynamic Young's modulus, $E^*$, a clamped material at a given temperature is made to vibrate at a particular amplitude and frequency. Thus, the temperature, vibration amplitude and frequency are con-trolled by the experimenter. A certain force is then measured in the experiment (the dependent variable), and using the appropriate beam theory one can find an inverse model for $E^*$. Experiments were done at the J. J. Pickle Research Campus of the University of Texas at Austin using the Q800 DMA instrument from TA Instruments. The experimental setup purposefully resembles cantilever beams. Indeed, Figure 4.2 shows a material sample in cantilever, clamped at both ends, where one clamp is static while the other clamp is free to move and vibrate at a given amplitude and frequency. It is at this moving clamp that the force is measured.

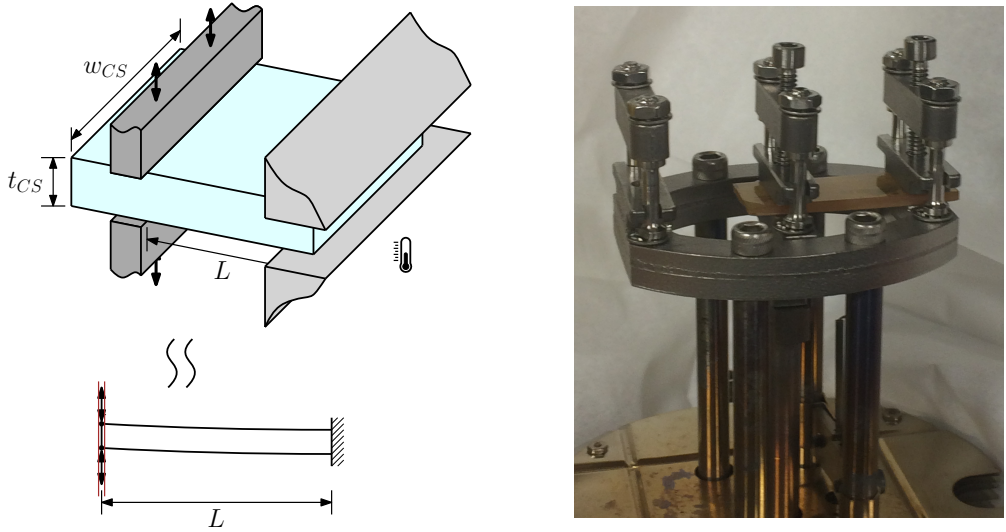

Figure 4.2: The single cantilever DMA experimental setup. The external clamp is statically fixed, while the central clamp, where a force is measured, moves vertically with a known amplitude and frequency. This whole setup lies inside a closed oven that carefully controls the temperature.

An inverse model for $E^*$ can be derived using Timoshenko beam theory. Consider a static linear elastic beam clamped at one end and with a point force applied at the other end, where additionally the cross-section remains parallel to the force (see Figure 4.2). This last condition represents the moving clamp where the force is being measured. Hence, this is *not* a typical cantilever beam (where one of the ends is free), but for simplicity it is still referred as such. Using Timoshenko beam theory [65, 66], the vertical displacement can be determined using the zero-angle boundary conditions at *both* ends and a zero-displacement in the clamped end. The resulting maximum displacement occurs where the force is applied and takes the value,

$$u_{\max} = \frac{FL^3}{12EI} + \frac{FL}{\kappa A_{CS} G} = \frac{FL^3}{12EI}\left(1 + \frac{2}{\kappa}(1+\nu)\left(\frac{t_{CS}}{L}\right)^2\right), \tag{4.1}$$

where $u_{\max}$ is the maximum vertical displacement of the beam, $F$ is the force applied, $L$ is the length between the clamped end and where the force is applied; $E$ and $G = \frac{E}{2(1+\nu)}$ are the Young's and shear moduli of the linear elastic material while $\nu$ is its Poisson's ratio; $A_{CS} = w_{CS} t_{CS}$, $w_{CS}$ and $t_{CS}$ are the cross-sectional area, width and thickness respectively; $I = \frac{w_{CS} t_{CS}^3}{12}$ is the second moment of area of the rectangular cross-section, and $\kappa$ is the Timoshenko shear coefficient. This equation obeys a correspondence principle with the time-harmonic equations of linear viscoelasticity [46], so that an inverse model of the form,

$$E^* = \frac{1}{\alpha_c} \frac{F^*_{exp}}{u^*_{max}} \frac{L^3}{\beta_c I} \left(1 + \frac{12}{5}(1+\nu^*)\left(\frac{t_{CS}}{L}\right)^2\right),$$

$$\alpha_c = 0.7616 - 0.02713\sqrt{\frac{L}{t_{CS}}} + 0.1083\ln\left(\frac{L}{t_{CS}}\right),$$

(4.2)

is utilized, where $\beta_c = 12$ in this single cantilever setting, and $\alpha_c$ is a correction factor accounting for 3D clamping effects, which is given by the manufacturer. For a rectangular cross-section, the Timoshenko shear coefficient is taken from the literature as $\kappa = \frac{5}{6}$ [42]. Here, $E^*$ is the dynamic Young's modulus, and note that both the experimental force and vibration amplitude, $F^*_{exp}$ and $u^*_{max}$, are now complex numbers. Note that $\frac{F^*_{exp}}{u^*_{max}} = \left|\frac{F^*_{exp}}{u^*_{max}}\right|e^{i\delta_{ph}}$, where $\delta_{ph}$ is an angle that represents the phase change between the oscillations of the force and the driving mechanical vibrations of the displacement. The values of temperature, vibration frequency, $|u^*_{max}|$, $|F^*_{exp}|\cos(\delta_{ph})$ and $\tan(\delta_{ph})$ are reported by the instrument. The distance $L$ here is the distance between the clamps themselves, *not* the distance between the midpoints of the clamps. The only limitation with this inverse model is that it assumes that the dynamic Poisson's ratio, $\nu^*$, is known. The ideal scenario is that either $\nu^*$ or the dynamic shear modulus, $G^*$, are known from a separate preceding experiment. In the latter case, where $G^*$ is known, note that $\nu* = \frac{E^*}{2G^*} - 1$, so an analogous expression for $E^*$ only in terms of $G^*$ can easily be derived from (4.2). If neither $\nu^*$ nor $G^*$ are experimentally known, it is usually assumed that $G^*$ has the same phase as $E^*$, so that $\nu^*$ is real-valued, and then an educated guess is made for $\nu^* \in \mathbb{R}$.

There is a second experimental setup which involves the same instrument, but with the beam arranged in a double cantilever, with two external static clamps at both ends and a middle moving clamp. This can be seen in Figure 4.3. The inverse model is actually the same as that given in (4.2), but with $\beta_c = 24$ in the double cantilever setting, and where the distance $L$ is the same as in the single cantilever case since it represents the distance between the edge of the external clamp and closest edge of the middle clamp (as seen in Figure 4.3). Hence, the actual distance over which the material is being deformed is $L_d = 2L$.

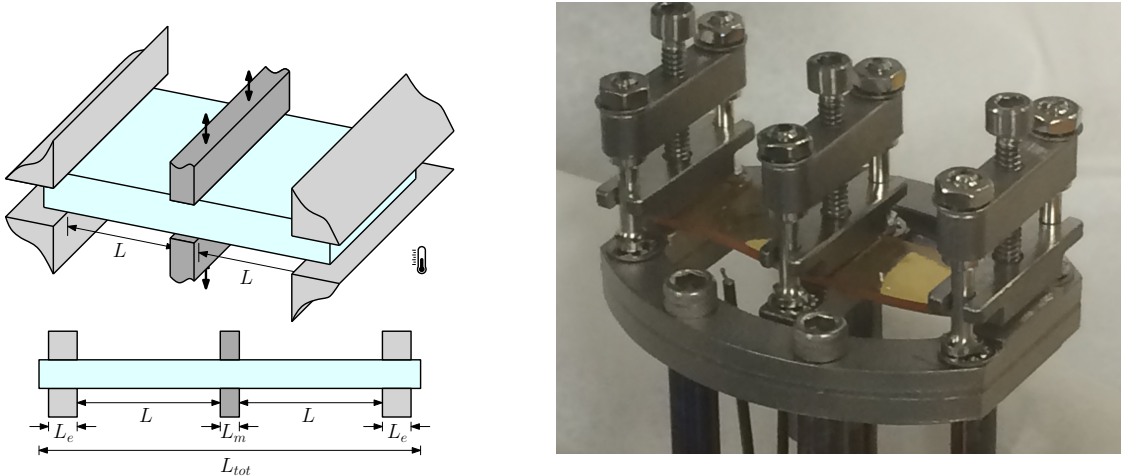

Figure 4.3: The double cantilever DMA experimental setup. The two external clamps are statically fixed, while a force is measured at the central clamp which moves at a controlled amplitude and frequency.

For the sake of brevity, results of only one example of each setup will be shown here. In the single cantilever case, silicone at 30.0°C was tested at 4Hz with an amplitude of vibration of $|u^*_{max}| = 15\mu m$, where the relevant part of the sample had dimensions of $L = 17.5$mm, $w_{CS} = 11.8$mm and $t_{CS} = 1.63$mm. The measured force from the experiment was $|F^*_{exp}|\cos(\delta_{ph}) = 0.1064$N with $\tan(\delta_{ph}) = 0.0384$. In the double cantilever case, epoxy at 22.4°C was tested at 40Hz with an amplitude of vibration of $|u^*_{\max}| = 15\mu$m, where the relevant part of the sample had dimensions of $L_d = 2L = 35.0$mm, $w_{CS} = 13.2$mm and $t_{CS} = 2.05$mm. The measured force from the experiment was $|F^*_{exp}|\cos(\delta_{ph}) = 0.7248$N with $\tan(\delta_{ph}) = 0.00869$. In both experiments it was assumed that $\nu^* = 0.33$ (see [30], but higher values are also found in [63]), so using the inverse model in (4.2) with $\beta_c = 12$ and $\beta_c = 24$ respectively, it was possible to calculate $E^*$.

Next, the dynamic stiffness tensor, $\mathsf{C}^*$, was computed using the values of $E^*$ and $\nu^*$, and the experiments were then simulated computationally. Here, it is important to mention that the the middle clamp measures $L_m = 6.35$mm, while the two external clamps measure $L_e = 7.625$mm each, as observed from Figure 4.3. Thus the samples themselves (both in the experiment and the simulated geometry) are typically longer than $L_e + L + L_m = 31.475$mm in the single cantilever case and $2L_e + 2L + L_m = 56.6$mm in the double cantilever case. The samples used for the numerical results were 40mm for the single cantilever and $L_{tot} = 60$mm for the double cantilever. The densities of the silicone and epoxy resins were assumed to be 1134kgm$^{-3}$ and 1250kgm$^{-3}$ respectively. The force, which is the quantity of interest, was calculated a posteriori by integrating the vertical traction, $(\hat{\sigma}_{n,h})_3$, over the area where the moving clamp made contact with the sample. The numerically computed force, $F^*_h$, was then compared with the actual measured force from the experiment, $F^*_{exp}$. The results for different values of $p$ and with $\Delta p = 1$ are shown in Figure 4.4.

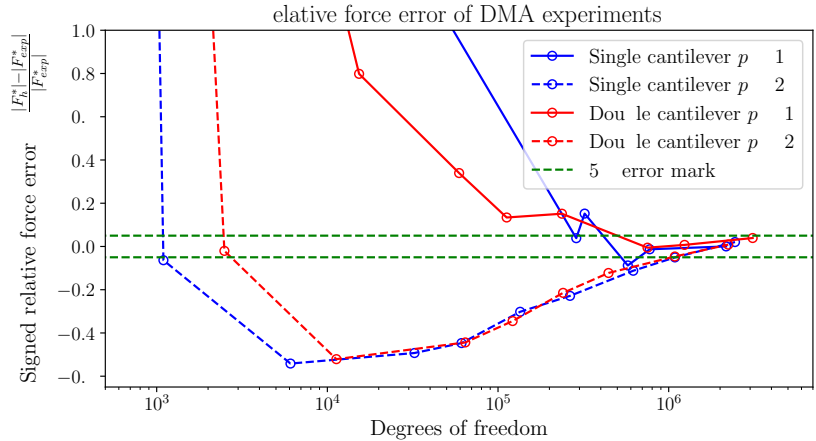

Figure 4.4: Convergence of the magnitude of the computed force, $F^*_h$, to the real experimental value measured from DMA experiments on different setups, $F^*_{exp}$. The single cantilever results correspond to a silicone sample, while those of the double cantilever correspond to an epoxy sample.

The magnitude of the force appears to converge to within 5% of the experimental value with both the single and double cantilever setups. This is as good as one can hope for from the validation point of view, and it confirms that the equations do indeed model the actual physical behavior observed experimentally. These results seem to suggest that the value of $p = 2$ does not offer a significant advantage over $p = 1$ to obtain the desired outcome, but further research on this matter might be necessary, as a different quantity of interest might produce very different results. With respect to the phase error in $\tan(\delta_{ph})$, the simulations show virtually no error even from the first computation. This is probably due to the assumption that $\nu^* \in \mathbb{R}$ is real-valued, but perhaps a less trivial convergence behavior would be observed if this hypothesis were to be dropped.

The results in Figure 4.4 were obtained with adaptivity driven by the arbitrary-$p$ residual-based a posteriori error estimator described in (3.13), which is innate to the DPG methodology. Otherwise, it would have been prohibitively expensive to obtain the same results via uniform refinements. Indeed, from the physics of the problem, it is intuitive to notice that most of the stress will be concentrated in the areas close to where the clamps are holding the material. The computations confirm this, as can be observed from Figure 4.5, where it is clear that not only the stress is localized there, but that the adaptivity scheme is refining in precisely that area, which is where the force will be computed from. Thus, adaptivity is fundamental for this problem which has localized solution features, and this justifies to a degree the use of the numerical method proposed here.

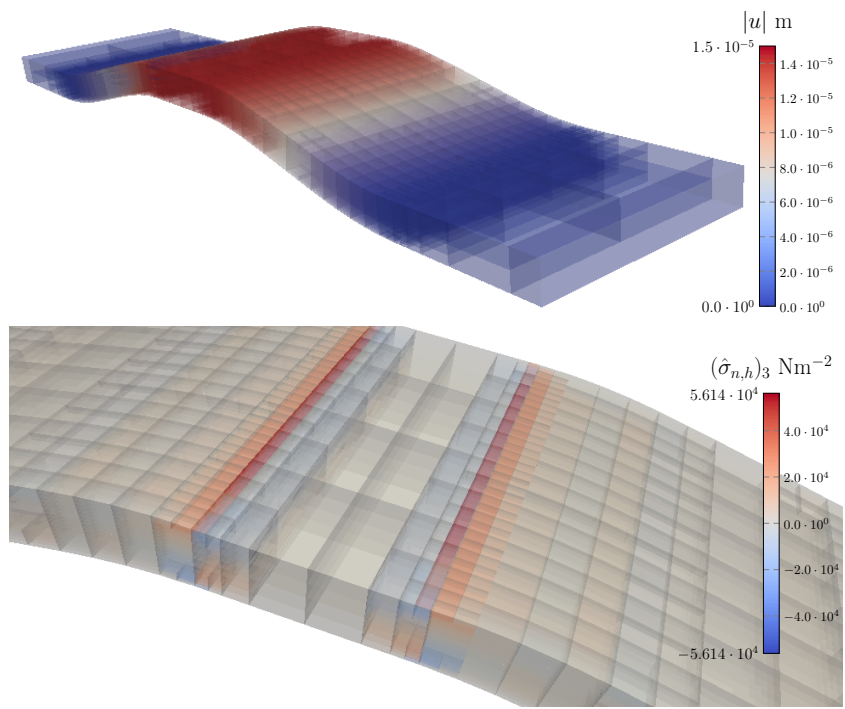

Figure 4.5: Numerical results with the double cantilever setup with $p = 1$ and after 4 isotropic adaptive refinements. The displacement is warped by a factor of 4000 for clarity. The vertical traction seems to be concentrated at the edges of the middle clamp, and adaptive refinements do seem to focus on that area.

# 5 Conclusions

A DPG finite element method was implemented for the time-harmonic equations of linear viscoelasticity. The method discretizes a broken primal variational formulation of the equation, which was proved to be well-posed in the infinite-dimensional setting. As part of this proof, the well-posedness of the classical primal variational formulation of linear viscoelasticity was also rigorously established. Moreover, the numerical method itself was shown to be stable and convergent under certain conditions, and this included a full $hp$-convergence analysis. A completely natural a posteriori error estimator for arbitrary-$p$ which is used to drive adaptivity is also included as part of the method. The method was verified using a smooth manufactured solution, where the expected $h$-convergence rates of the form $h^p$ where corroborated for various values of $p$. Moreover, the verification tests displayed exponential $p$-convergence of the form $\exp(-b\,p^{1.25})$.

Additionally, DMA experiments to determine the dynamic Young's modulus, $E^*$, were performed on different materials and with distinct experimental setups: single and double cantilever. The computational results validated the calibration model to within 5% error of the quantity of interest. Moreover, the simulated stress was very concentrated on certain parts of the domain, so having a good adaptivity scheme was crucial to obtain the desired result. In this sense, the numerical method was extremely convenient, since it already came with its own a posteriori error estimator.

Looking forward, more complicated validation studies could be tackled, where the quantities of interest may vary in nature. The built-in a posteriori error estimator is designed to drive down the residual, but may not be optimal in accelerating the convergence of a particular quantity of interest. In this sense, this could lead to investigating goal-driven adaptivity schemes within the context of the DPG methodology. When the linear system size is large, computations may become prohibitive, so it would be useful to make improvements to reduce the system size as much as possible and to support parallel computing within the solvers. Finally, for more interesting cases closer to the glass transition temperature of the materials in question, the results from the computations might improve if the actual value of the dynamic Poisson's ratio, $\nu^*$, or the dynamic shear modulus, $G^*$, are used in the calibration inverse model, but this requires separate DMA experiments to be completed, which might be a future endeavor.

**Acknowledgements.** This work was partially supported with grants by ONR (N00014-15-1-2496), NSF (DMS-1418822), and AFOSR (FA9550-12-1-0484).

[40] Houston, P., Schwab, C., and Süli, E. (2002). Discontinuous *hp*-finite element methods for advection-diffusion-reaction problems. *SIAM J. Numer. Anal.*, 39(6):2133–2163.

[41] Ihlenburg, F. (1998). *Finite Element Analysis of Acoustic Scattering*, volume 132 of *Applied Mathematical Sciences*. Springer-Verlag, New York.

[42] Kaneko, T. (1975). On Timoshenko's correction for shear in vibrating beams. *J. Phys. D: Appl. Phys.*, 8(16):1927–1936.

[43] Keith, B., Fuentes, F., and Demkowicz, L. (2016). The DPG methodology applied to different variational formulations of linear elasticity. *Comput. Methods Appl. Mech. Engrg.*, 309:579–609.

[44] Keith, B., Knechtges, P., Roberts, N. V., Elgeti, S., Behr, M., and Demkowicz, L. (2017a). An ultraweak DPG method for viscoelastic fluids. *J. Nonnewton. Fluid. Mech.*

[45] Keith, B., Petrides, S., Fuentes, F., and Demkowicz, L. (2017b). Discrete least-squares finite element methods. *ArXiv e-prints*, arXiv:1705.02078 [math.NA].

[46] Lakes, R. (2009). *Viscoelastic Materials*. Cambridge University Press, Cambridge.

[47] Lakes, R. and Wojciechowski, K. W. (2008). Negative compressibility, negative Poisson's ratio, and stability. *Phys. Status Solidi B*, 245(3):545–551.

[48] McLean, W. (2000). *Strongly Elliptic Systems and Boundary Integral Equations*. Cambridge University Press, Cambridge.

[49] Muga, I. and van der Zee, K. G. (2015). Discretization of linear problems in Banach spaces: residual minimization, nonlinear Petrov–Galerkin, and monotone mixed methods. *ArXiv e-prints*, arXiv:1511.04400 [math.NA].

[50] Nagaraj, S., Petrides, S., and Demkowicz, L. (2017). Construction of DPG Fortin operators for second order problems. *Comput. Math. Appl.*

[51] Nédélec, J. C. (1980). Mixed finite elements in $\mathbb{R}^3$. *Numer. Math.*, 35:315–341.

[52] Neff, P., Pauly, D., and Witsch, K.-J. (2015). Poincaré meets Korn via Maxwell: Extending Korn's first inequality to incompatible tensor fields. *J. Differential Equations*, 258(4):1267–1302.

[53] Niemi, A. H., Collier, N., and Calo, V. M. (2013). Automatically stable discontinuous Petrov-Galerkin methods for stationary transport problems: Quasi–optimal test space norm. *Comput. Math. Appl.*, 66(10):2096–2113.

[54] Oden, J. T. (1972). *Finite Elements of Nonlinear Continua*. McGraw-Hill, New York. Reprinted in 2006 by Dover Publications, Inc., Mineola, NY.

[55] Oden, J. T. and Armstrong, W. H. (1971). Analysis of nonlinear, dynamic coupled thermoviscoelasticity problems by the finite element method. *Comput. Struct.*, 1(4):603–621.

[56] Oden, J. T. and Demkowicz, L. (2010). *Applied Functional Analysis*. Chapman & Hall/CRC Press, New York, 2nd edition.